\documentclass{article}

\usepackage[accepted]{icml2020}
\usepackage{microtype}
\usepackage{fancyhdr} 
\usepackage{mathrsfs}
\usepackage{graphicx}
\usepackage{subfigure}
\usepackage{amssymb}
\usepackage{booktabs} 
\usepackage{bm}
\usepackage{url}            
\usepackage{booktabs}       
\usepackage{amsfonts}       
\usepackage{nicefrac}       
\usepackage{microtype}      
\usepackage{bm}
\usepackage{amsmath}
\usepackage{amsfonts}
\usepackage{amssymb}
\usepackage{mathtools}
\usepackage{titlesec}
\usepackage{wrapfig}
\usepackage{caption}
\usepackage{natbib}
\usepackage{enumitem}
\usepackage{setspace}
\usepackage{caption}

\usepackage[colorlinks=true, bookmarksopen,
            pdfsubject={algorithms},
            linkcolor={blue},
            anchorcolor={black},
            citecolor={black},
            filecolor={magenta},
            menucolor={black},
            plainpages=false,pdfpagelabels,
            urlcolor={blue}]{hyperref}

\newtheorem{thm}{Theorem}

\newtheorem{rem}{Remark}

\icmltitlerunning{TransNet: Transferable Neural Networks for PDEs}

\begin{document}

\twocolumn[
\icmltitle{TransNet: Transferable Neural Networks for Partial Differential Equations}



\icmlsetsymbol{equal}{*}

\begin{icmlauthorlist}
\icmlauthor{Zezhong Zhang}{FSU}
\icmlauthor{Feng Bao}{FSU}
\icmlauthor{Lili Ju}{USC}
\icmlauthor{Guannan Zhang}{CSMD,equal}
\end{icmlauthorlist}

\icmlaffiliation{CSMD}{Computer Science and Mathematics Division, Oak Ridge National Laboratory, TN 37831, USA}
\icmlaffiliation{FSU}{Department of Mathematics, Florida State University, Tallahassee, FL 32306, USA}
\icmlaffiliation{USC}{Department of Mathematics, University of South Carolina, Columbia, SC 29208, USA}

\icmlcorrespondingauthor{Guannan Zhang}{zhangg@ornl.gov}

\icmlkeywords{Machine Learning, ICML}
\vskip 0.3in
]



\printAffiliationsAndNotice{\hspace{-0.4cm} * Corresponding author\\} 

\begin{abstract}
Transfer learning for partial differential equations (PDEs) is to develop a pre-trained neural network that can be used to solve a wide class of PDEs. Existing transfer learning approaches require much information of the target PDEs such as its formulation and/or data of its solution for pre-training. In this work, we propose to construct transferable neural feature spaces from purely function approximation perspectives without using PDE information. The construction of the feature space involves re-parameterization of the hidden neurons and uses auxiliary functions to tune the resulting feature space. Theoretical analysis shows the high quality of the produced feature space, i.e., uniformly distributed neurons. Extensive numerical experiments verify the outstanding performance of our method, including significantly improved transferability, e.g., using the same feature space for various PDEs with different domains and boundary conditions, and the superior accuracy, e.g., several orders of magnitude smaller mean squared error than the state of the art methods.  
\end{abstract}

\section{Introduction}
\label{sec:intro}
Rapid advancement of deep learning has attracted significant attention of researchers to explore how to use deep learning to solve scientific and engineering problems. Since numerical solutions of partial differential equations (PDEs) sits at the heart of many scientific areas, there is a surge of studies on how to use neural networks to leverage data and physical knowledge to solve PDEs \citep{raissi2019physics,e2018ritz,long2018pde,zang2019weak,li2020fourier,li2020multipole,lu2021learning,gin2021deepgreen,zhang2021ifnn,Teng2022gfnet,CLARKDILEONI2023111793}. The neural network-based methods have several advantages over traditional numerical methods (e.g., finite element, finite difference and finite volume), such as avoiding the need for numerical integration, generating differentiable solutions, exploiting advanced computing capabilities, e.g., GPUs. Nevertheless, a major drawback of these deep learning methods for solving PDEs is high computational cost associated with the neural network training/retraining using stochastic gradient descent (SGD). One of the popular strategies to alleviate this issue is transfer learning. 

Transfer learning for PDEs is to develop a pre-trained neural network that can be effectively re-used to solve a PDE with multiple coefficients or in various domains, or to solve multiple types of  PDEs. When transfer a pre-trained neural network from one scenario to another, the feature space, e.g., the hidden layers, are often frozen or slightly perturbed, which can dramatically reduce the training overhead by orders of magnitude. However, existing transfer learning approaches for PDEs, e.g., \citep{lu2021learning, li2020fourier,Chakraborty2020TransferLB, Desai2021OneShotTL}, require information/knowledge of the target family of PDEs to pre-train a neural network model. The needed information could be the analytical definitions of the PDEs including initial and boundary conditions, and/or measurement data of the PDE's solution. These requirement not only leads to time-consuming simulation data generation using other PDE solvers, but also limits the transferability of the pre-trained neural network (i.e., the pre-trained network is only transferable to the same or similar type of PDEs that are used for pre-training). 

To overcome the above challenges, in this paper we propose a transferable neural network (TransNet) to improve the transferability of neural networks for solving PDEs. The key idea is construct a pre-trained neural feature space without using any PDE information, so that the pre-trained feature space could be transferred to a variety of PDEs with different domains and boundary conditions. We limit our attention to single-hidden-layer fully-connected neural networks, which have sufficient expressive power for low-dimensional PDEs that are commonly used in science and engineering fields. Specifically, we treat each hidden neuron as a basis function and re-parameterize all the neurons to separate the parameters that determine the neuron's location and the ones that control the shape (i.e., the slope) of the activation function. Then, we develop a simple, yet very effective, approach to generate uniformly distributed neurons in the unit ball, and rigorously prove the uniform neuron distribution. Then, the shape parameters of the neurons are tuned using auxiliary functions, i.e., realizations of a Gaussian process. The entire feature space construction (determining the hidden neurons' parameters) does not require the PDE's formulation or data of the PDE's solution. When applying the constructed feature space to a PDE problem, we only need to solve for the parameters of the output layer by minimizing the standard PDE residual loss. This can be done by either solving a simple least squares problem for linear PDE or combining a least squares solver with a nonlinear iterative solver, e.g., Pichard iteration, for nonlinear PDEs.

The major contributions of this work are summarized as
\vspace{-0.25cm}
\begin{itemize}[leftmargin=10pt]\itemsep-0.0cm
    \item We develop transferable neural feature spaces that are 
    {independent} of any PDE, and can be applied to effectively solve various linear and nonlinear PDE problems. 

    \item We theoretically and computationally prove the uniform distribution of the hidden neurons, viewed as global non-orthogonal basis, for the proposed TransNet in the unit ball of any dimension. 

    \item We demonstrate the superior accuracy and efficiency of the proposed TransNet for solving PDEs, e.g., the mean square errors of TransNet are several orders of magnitudes smaller than those by the state-of-the-art methods. 
    
\end{itemize}

\section{Related work}
Studies on using neural networks for solving PDEs can be traced back to some early works, e.g.,  \citep{dissanayake1994neural,lagaris1998artificial}. Recent advances mostly have been focused on physics-informed neural network (PINN). The general idea of PINN is to represent the PDE's solution by a neural network, and then train the network by minimizing certain measurement of the PDE's residual at a set of samples in the domain of computation. Several improvements on the training and sampling were proposed in \citep{lu2021deepxde,anitescu2019artificial,zhao2020solving,krishnapriyan2021characterizing}. Besides direct minimizing the PDE's residual, there are studies on how to combine traditional PDE solvers with neural networks. For example, the deep Ritz method \citep{e2018ritz} uses the variational form of PDEs and combines the stochastic gradient descent with numerical integration to train the network; the deep Galerkin method \citep{sirignano2018} combines the Galerkin method with machine learning; 
the PDE-Net \citep{long2018pde,long2019pde} uses a stack of neural networks to approximate the PDE solutions over a multiple of time steps.

Another type of deep learning method for PDEs is to use neural networks to learn a family of PDE operators, instead of a single equation. The Fourier neural operator (FNO)  \citep{li2020fourier} parameterizes the integral kernel in Fourier space and is generalizable to different spatial/time resolutions. The DeepONet \citep{lu2021learning} extends the universal approximation theorem \citep{chen1995universal} to deep neural networks, and its variant \citep{wang2021learning} further reduces the amount of data needed for training. The physics-informed neural operator (PINO) \citep{li2021physics} combines operator learning with function approximation to achieve higher accuracy.
MIONet \cite{JML2022MIONet} was proposed to learn multiple-input operators via tensor product basd on low-rank approximation.

Random feature models have also been used to solve PDEs \citep{Sun2018OnTA,https://doi.org/10.48550/arxiv.2212.05591} or learn PDE operators \citep{doi:10.1137/20M133957X}. The theory of random feature models for function approximation was developed due to its natural connection with kernel methods \citep{9495136,10.5555/3122009.3122030}. The proposed TransNet can be viewed as an improved random feature model for PDEs from two perspectives: (1) the re-parameterization of the hidden neurons to separate the parameters that determine locations of the neurons and the ones that control the activation function slope, (2) the usage of auxiliary functions to tune the neural feature space, which makes a critical contribution to the improvement of the accuracy of TransNet in solving PDEs.

\section{Transferable neural networks for PDEs}\label{sec:method}

\subsection{Problem setting and background}\label{sec:backgroud}
We introduce the problem setup for using neural networks to solve partial differential equations. The PDE of interest can be presented in a general formulation, i.e., 
\begin{equation}\label{eq:pde}
\left\{
\begin{aligned}
 & \mathcal{L}(u(\bm y)) = f(\bm y)\;\;  \text{ for }\; \bm y \in \Omega,\\
 & \mathcal{B}(u(\bm y)) =  g(\bm y)\;\;  \text{ for }\; \bm y \in \partial \Omega,
\end{aligned}
\right.
\end{equation}
where $\Omega \subset \mathbb{R}^{d}$ with the boundary $\partial \Omega$ is the spatial-temporal bounded domain under consideration, 
$\bm y := (\bm x, t) = (x_1, \ldots, x_{d-1}, t)^{\top}$ is a column vector includes both spatial and temporal variables, $u$ denotes the unknown solution of the PDE, $\mathcal{L}(\cdot)$ is a differential operator, $\mathcal{B}(\cdot)$ is the operator defining the initial and/or boundary conditions, $ f(\bm y)$ and $g(\bm y)$ are the right hand sides associated with the operators $\mathcal{L}(\cdot)$ and $\mathcal{B}(\cdot)$, respectively. For notational simplicity, we assume that the solution is a scalar function; the proposed method can be extended to vector-valued functions without any essential difficulty. 
We limit our attention to the single-hidden-layer fully-connected neural networks, denoted by
\begin{equation}\label{eq:fc}
    u_{\rm NN}(\bm y) := \sum_{m = 1}^M \alpha_m\, \sigma( \bm w_m \bm y + b_m) + \alpha_0,
\end{equation}
where $M$ is the number of hidden neurons, the row vector $\bm w_m = (w_{m,1}, \ldots, w_{m,d})$ and the scalar $b_m$ are the weights and bias of the $m$-th hidden neuron, the row vector $\bm \alpha = (\alpha_0, \alpha_1, \ldots, \alpha_M)$ includes the weights and bias of the output layer, and $\sigma(\cdot)$ is the activation function. As demonstrated in Section \ref{sec:exp}, this type of neural networks have sufficient expressive power for solving a variety of PDEs with satisfactory accuracy.

A typical method \citep{2021NatRP...3..422K} for solving the PDE in Eq.~\eqref{eq:pde} is to directly parameterize the solution $u(\bm y)$ as a neural network $u_{\rm NN}(\bm y)$ in Eq.~\eqref{eq:fc} and optimize the neural network's parameters by minimizing the PDE residual loss, e.g., $L(\bm y) = \|\mathcal{L}(u(\bm y)) - \mathcal{L}(u_{\rm NN}(\bm y))\|_2 + \|\mathcal{B}(u(\bm y)) - \mathcal{B}(u_{\rm NN}(\bm y)) \|_2$, at a set of spatial-temporal locations. Despite the good performance of these approaches in solving  PDE problems, its main drawback is the {\em limited transferability} because of the high computational cost of gradient-based re-training and hyperparameter re-tuning. When there is any change to the operators $\mathcal{L}(\cdot), \mathcal{B}(\cdot)$, the right-hand-side functions $f(\bm y), g(\bm y)$, or the shape of the domain $\Omega$, the neural network $u_{\rm NN}(\bm y)$ often needs to be re-trained using gradient-based optimization (even though the current parameter values could provide a good initial guess for the re-training), or the hyperparameters associated with the network and the optimizer need to be re-tuned. In comparison, the random feature models require much lower re-training cost, which has been exploited in learning operators \citep{doi:10.1137/20M133957X} and dynamical systems \citep{NEURIPS2021_72fe6f9f, https://doi.org/10.48550/arxiv.2212.05591}.


\subsection{The neural feature space} We can treat each hidden neuron $\sigma(\bm w_m \bm y + b_m)$ as a nonlinear feature map from the space of $\bm y \in \mathbb{R}^{d}$ to the output space $\mathbb{R}$. From the perspective of approximation theory, the set of hidden neurons $\{\sigma(\bm w_m \bm y + b_m)\}_{m=1}^M$ can be viewed as a globally supported basis in $\mathbb{R}^d$. The {\em neural feature space}, denoted by $\mathcal{P}_{\rm NN}$, can be defined by the linear space expanded by the basis $\{\sigma(\bm w_m \bm y + b_m)\}_{m=1}^M$, i.e., 
\begin{equation}\label{eq:fs}
    \mathcal{P}_{\rm NN} = {span}\Big\{1, \sigma(\bm w_1 \bm y + b_1), \ldots, \sigma(\bm w_M \bm y + b_M) \Big\},
\end{equation}
where the constant basis corresponds to the bias of the output layer. Then, the 
neural network in Eq.~\eqref{eq:fc} lives in the linear space, i.e., $
u_{\rm NN}(\bm y) \in \mathcal{P}_{\rm NN}.$
In other words, the neural network approximation can be viewed as a spectral method with {\em non-orthogonal} basis, and the parameters $\bm \alpha$ in Eq.~\eqref{eq:fc} of the output layer of $u_{\rm NN}(\bm y)$ contains the coefficients of the expansion in the neural feature space $\mathcal{P}_{\rm NN}$. 

In the PINN methods, the neural feature space $\mathcal{P}_{\rm NN}$ and the coefficient $\bm \alpha$ are trained simultaneously using stochastic gradient descent methods, which often leads to a non-convex and ill-conditioned optimization problem. It has been shown that the non-convexity and ill-conditioning in the neural network training are major reasons of unsatisfactory accuracy of the trained neural network. A natural idea to reduce the complexity of the training is to decouple the training of $\mathcal{P}_{\rm NN}$ from that of $\bm \alpha$. For example, in random feature models, $\mathcal{P}_{\rm NN}$ is defined by randomly generating the parameters$\{(\bm w_m, b_m)\}_{m=1}^M$ from a user-defined probability distribution; the coefficients $\bm \alpha$ can then be obtained by solving a linear system when the operators $\mathcal{L}$, $\mathcal{B}$ in Eq.~\eqref{eq:pde} are linear.  However, the numerical experiments in Section \ref{sec:exp} show that the random feature model based on Eq.~\eqref{eq:fc} converges very slowly with the increase of the number of features. This drawback motivates us to develop a methodology to customize the neural feature space $\mathcal{P}_{\rm NN}$ to improve the accuracy, efficiency and transferability of $u_{\rm NN}$ in solving PDEs.

\subsection{Constructing the transferable neural feature space}\label{sec:par}
%
This section contains the key ingredients of the proposed TransNet. The goal is to construct a single neural feature space $\mathcal{P}_{\rm NN}$ that can be used to solve various PDEs in different domains.

\subsubsection{Re-parameterization of $\mathcal{P}_{\rm NN}$}\label{sec:repar}
The first step is to re-parameterize the hidden neuron $\sigma(\bm w_m \bm y + b_m)$, viewed as a basis function in $\Omega$, to separate the components that determine the {\em location} of the neuron and the components that control the {\em shape} of the neuron. 

The idea of handling the locations of the basis functions is inspired by the studies on activation patterns of ReLU networks. When $\sigma$ is the ReLU function, there is a {\em partition hyperplane} defined by 
\begin{equation}\label{eq:plane}
    w_{m,1}  y_1 + w_{m,2} y_2 + \cdots + w_{m,d} y_d + b_m = 0
\end{equation}
that separates the activated and inactivated regions for this neuron. The intersections of multiple partition hyperplanes associated with different neurons define a linear region of ReLU network. Studies have shown that the expressive power of a ReLU network is determined by the number of linear regions and the distribution of those linear regions. In principle, the more {\em uniformly distributed} linear regions in the domain $\Omega$, the more expressive power the ReLU network has. For other activation functions, e.g., $tanh(\cdot)$ that is widely used in solving PDEs due to its smoothness, the partition hyperplane in Eq.~\eqref{eq:plane} can be used to describe the geometric property of the neuron.

\begin{figure*}[t]
    \centering
    \includegraphics[width=\textwidth]{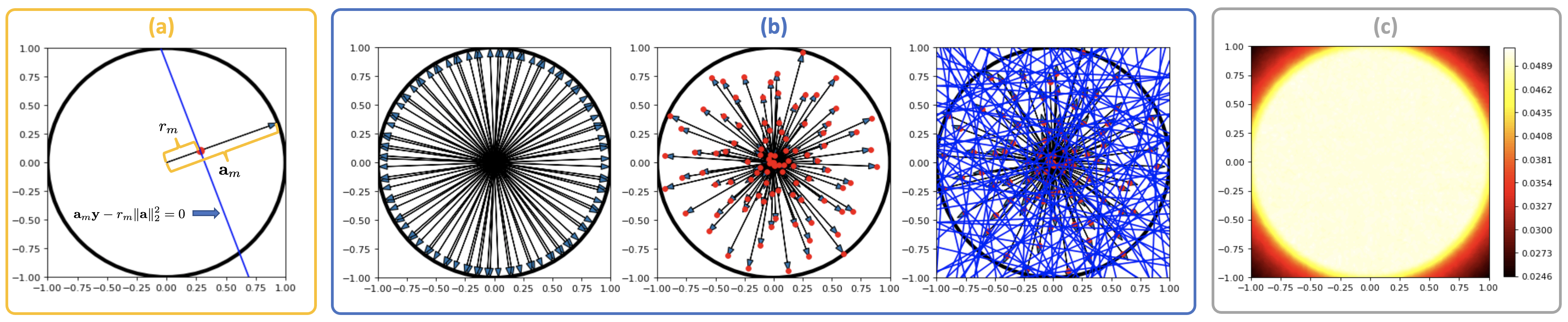}
    \vspace{-0.7cm}
    \caption{{\bf (a)} Illustrates how the re-parameterization in Eq.~\eqref{eq:plane2} characterizes the location of a neuron. The blue line is the plane where $\tanh(\cdot) = 0$, 
    $\bm a_m$ (the arrow) is the normal direction of the plane, the red dot is the location $r_m \bm a_m$ that the plane passes, $r_m$ is the distance between the origin and the plane. {\bf (b)} illustrates how to generate uniformly distributed neurons in the unit ball. The first step in {\bf (b)-left} is to generate the normal directions $\{\bm a_m\}_{m=1}^M$ uniformly distributed on unit sphere; the second step in {\bf (b)-middle} is to generated $\{r_m\}_{m=1}^M$ uniformly from $[0,1]$ defining the locations the neurons' partition hyperplanes will pass; the blue lines in {\bf (b)-right} show the distribution of the partition hyperplanes. {\bf (c)} shows the density function $D_M(\bm y)$ with $\tau = 0.05$ in Eq.~\eqref{eq:den} for a set of neurons generated using our approach. We can see that our approach provides a uniformly distributed neurons in the ball $B_{1-\tau}(\bm 0)$, which is consistent with Theorem \ref{thm1}.} 
    \label{fig1}
\end{figure*}

Specifically, let us re-write Eq.~\eqref{eq:plane} into the following point-slope form:
\begin{equation}\label{eq:plane2}
    \gamma_m\big(a_{m,1} (y_1 - r_m a_{m,1}) + \cdots + a_{m,d} (y_d - r_m a_{m,d})\big) = 0,
\end{equation}
where $\bm a_m = (a_{m,1}, \ldots, a_{m,d})$ is a unit vector, i.e., $\|\bm a_m \|_2 = 1$, $r_m > 0$ and $\gamma_m \in \mathbb{R}$ are two scalar parameters for the $m$-th neuron.
We can relate Eq.~\eqref{eq:plane2} to Eq.~\eqref{eq:plane} by 
\begin{equation}\label{eq:rep}
\left\{\begin{split}
   w_{m,i} &= \gamma_m a_{m,i}, \;\;\;i=1,\cdots,d,\\
   b_m &= - \gamma_m \sum_{i=1}^d a_{m,i}^2 r_m,  
\end{split}\right.
\end{equation}
which shows the desired geometric properties of the partition hyperplane in Eq.~\eqref{eq:plane}. In terms of the location, the unit vector $\bm a_m$ is the normal direction of the partition hyperplane in $\mathbb{R}^{d}$, the vector $(r_m a_{m,1}, \ldots, r_m a_{m,d})$ indicates a point that the hyperplane passes, 
$r_m$ is the distance between the origin and the partition hyperplane. An illustration is shown in Figure \ref{fig1}{\color{blue}(a)}. In terms of the shape, the constant $\gamma_m$ determines the steepness of the slope of the activation function along the normal direction $\bm a_m$. Thus, the re-parameterization in Eq.~\eqref{eq:plane2} successfully separates the parameters determining location from the ones determining the shape.


\subsubsection{Generating uniformly distributed neurons for $\mathcal{P}_{\rm NN}$}\label{sec:uniform}
The second step of constructing $\mathcal{P}_{\rm NN}$ is to determine the parameters $\{(\bm a_m, r_m)\}_{m=1}^M$ in Eq.~\eqref{eq:plane2}, such that all the neurons are uniformly distributed in $\Omega$. 
We assume $\Omega$ is a {\em unit ball}, i.e., $B_{1}(\bm 0)=\{\bm y:\|\bm y\|_2 \le 1\} \subset \mathbb{R}^d$ in this subsection. To proceed, we need to define a density function that measures the neuron distribution. 
%
 %
 For a given $\bm y \in \Omega$, the distance between $\bm y$ and the partition hyperplane in Eq.~\eqref{eq:plane2} is given by  
\begin{equation}\label{eq:dist}
    dist(\bm y, m) = |\bm a_m (\bm y - r_m \bm a_m)|,
\end{equation}
for $m = 1, \ldots, M$. We use this distance to define how close the point $\bm y$ to the $m$-th neuron. 
The density function, denoted by $D_M(\bm y)$, is defined using the above distance, i.e., 
\begin{equation}\label{eq:den}
    D_M(\bm y) = \frac{1}{M}\sum_{m=1}^M \mathbf{1}_{dist(\bm y, m) < \tau} (\bm y),
\end{equation}
where $\mathbf{1}_{dist(\bm y, m) < \tau} (\bm y)$ is the indicator function of the event that the distance between $\bm y$ and the $m$-th neuron is smaller than a prescribed tolerance $\tau > 0$. Intuitively, $D_M(\bm y)$ measures the percentage of neurons whose partition hyperplane in Eq.~\eqref{eq:plane} intersect the ball (with radius $\tau$) around $\bm y$.

Next we propose the following approach, illustrated in Figure \ref{fig1}{\color{blue}(b)}, to generate the parameters $\{(\bm a_m, r_m)\}_{m=1}^M$. Specifically, we first generate the normal directions $\{\bm a_m\}_{m=1}^M$ uniformly distributed on the $d-1$-dimensional unit sphere. Note that when $d>2$, sampling uniformly in the angular space in the hyperspherical coordinate system does not lead to uniformly distributed samples on the unit sphere. This is known as the sphere point picking problem. To overcome this issue, we draw samples from the $d$-dimensional Gaussian distribution in the Cartesian coordinate system, and normalize the samples to unit vectors to obtain $\{\bm a_m\}_{m=1}^M$. Then, we generate $\{r_m\}_{m=1}^M$ uniformly from $[0,1]$ using the Monte Carlo method. The following theorem shows that our approach provides a set of uniformly distributed neurons in $\Omega$, where the density is measured by $D_M(\bm y)$ in Eq.~\eqref{eq:den}. 

\begin{thm}[Uniform neuron distribution]\label{thm1}
Given the re-parameterization in Eq.~\eqref{eq:plane2}, if $\{\bm a_m\}_{m=1}^M$ are uniformly distributed random vectors on the $d$-dimensional unit sphere, i.e., $\|\bm a_m\|_2 = 1$, and $\{r_m\}_{m=1}^M$ are uniformly distributed random variables in $[0,1]$, then, for a fixed $\tau \in (0,1)$,
\[
\mathbb{E}[D_M(\bm y)] = \tau\; \text{ for any}\; \|\bm y\|_2 \le 1-\tau,
\]
where $D_M(\bm y)$ is the density function defined in Eq.~\eqref{eq:den}.
\end{thm}

The proof is given in Appendix \ref{sec:proof}; an illustration of the density function is given in Figure \ref{fig1}{\color{blue}(c)}. This result is a little surprising that the distribution of $\{r_m\bm a_m\}_{m=1}^M$, i.e., the red dots in Figure \ref{fig1}{\color{blue}(b)-middle}, are not uniformly distributed in the ball $B_{1-\tau}(\bm 0)$, but the density function $D_M(\bm y)$ is a constant in the ball $B_{1-\tau}(\bm 0)$.

\begin{rem}[The dimentionality]\label{rem3}
Even though Theorem \ref{thm1} holds for any dimension $d$, the number of neurons required to cover a high-dimensional unit ball still could be intractable. 
On the other hand, the majority of PDEs commonly used in science and engineering are defined in low-dimensional domains, e.g., 3D spatial domain + 1D time domain. In this scenario, the proposed method is effective and easy to implement,
as demonstrated in Section \ref{sec:exp}.
\end{rem}

\subsubsection{Tuning the shape of the neurons in $\mathcal{P}_{\rm NN}$ using auxiliary functions }\label{sec:train}
The third step is to tune the shape parameters $\{\gamma_m\}_{m=1}^M$ in Eq.~\eqref{eq:plane2} that controls the slope of the activation function.
The experimental tests in Section \ref{sec:basis} show that the slope parameters play a critical role in determining the accuracy of the neural network approximator $u_{\rm NN}$. 
For simplicity, we assume the same shape parameter value for all neurons, i.e., 
$
\gamma = \gamma_m \text{ for } m = 1, \ldots, M.
$
Because we intend to construct a feature space $\mathcal{P}_{\rm NN}$ that
can be used in multiple scenarios, e.g., various PDEs with different domains and boundary conditions, we do not want to tune the shape parameter $\gamma$ using any information about a specific PDE. 

Our idea is to use auxiliary functions that have similar or more complicated spatial-temporal variation frequency as the PDE solution to tune $\gamma$. Specifically, we propose to use realizations of Gaussian processes to generate the auxiliary functions. The advantage of Gaussian process is that one can control the variation frequency of its realizations by adjusting the correlation length. Additionally, the Guassian process is independent of the coordinate system. 
Let us denote by ${G}(\bm y| \omega, \eta)$ the Gaussian process, where $\omega$ represents the abstract random variable and $\eta$ is the correlation length. Given a correlation length, we first generate a set of realizations of the Gaussian process, denoted by $\{G(\bm y | \omega_k, \eta)\}_{k = 1}^K$. For each realization, define the MSE loss as
\begin{equation}\label{eq:mse}
\begin{aligned}
& \text{MSE}(u_{\rm NN}(\bm y), G(\bm y|\omega_k, \eta)) \\
= & 
\frac{1}{J}\sum_{j=1}^J \left[\sum_{m=1}^M \alpha_m \sigma(\bm w_m \bm y_j + b_m) + \alpha_0 - G(\bm y_j|w_k, \eta) \right]^2,
\end{aligned}
\end{equation}
where the parameters $\{\bm w_m\}_{m=1}^M$ and $\{b_m\}_{m=1}^M$ are already determined using the strategy in Section \ref{sec:uniform} and Eq.~\eqref{eq:rep}, and $J$ denotes the number of sample points. Unlike standard neural network training, the optimal coefficient $\bm \alpha$ that minimizing the MSE loss can be efficiently achieved  by solving the least squares problem. Hence, the shape parameter $\gamma$ can be tuned by solving the following one-dimensional optimization problem
%
%
%
\begin{equation}\label{eq:ls2}
\min_{\gamma} \left\{\sum_{k=1}^K \min_{\bm \alpha}\left[\text{MSE}({u}_{\rm NN}(\bm y), G(\bm y| \omega_k, \eta)) \right]\right\},
\end{equation}
where for each candidate $\gamma$, we solve $K$ least squares problems to compute the total loss. 

\begin{rem}[The choice of the correlation length]
There are two strategies to choose the correlation length $\eta$. One is to use the prior knowledge about the PDE. For example, for the Naveier-Stokes equations with low Reynolds' number, we know the solution will not have very high-frequency oscillation. The other  is to use an over-killing correlation length to ensure that the feature space has sufficient expressive power to solve the target PDE.
\end{rem}

\subsection{Applying TransNet to linear and nonlinear PDEs}
Once the neural feature space $\mathcal{P}_{\rm NN}$ is constructed and tuned, we can readily use it to solve PDE problems. Even though $\mathcal{P}_{\rm NN}$ is defined on the unit ball, i.e., $B_{1}(\bm 0)$, we can always place the (bounded) domain $\Omega$ for the target PDE in $B_{1}(\bm 0)$ by simple translation and dilation. Thus, the feature space can be used to handle PDEs defined in various domains, as demonstrated in Section \ref{sec:exp}.

{\bf Linear PDEs.}\; 
When $\mathcal{L}$ and $\mathcal{B}$ in Eq.~\eqref{eq:pde} are linear operators, the unknown parameters $\bm \alpha = (\alpha_0, \ldots, \alpha_M)$ in Eq.~\eqref{eq:fc} can be easily determined by solving the following least squares problem, i.e., 
\begin{equation}\label{eq:ls}\small
\begin{aligned}
& \min_{\bm \alpha}\Bigg\{ \frac{1}{J_1}\sum_{j=1}^{J_1}\left[\sum_{m=1}^M\alpha_m\, \mathcal{L}(\sigma( \bm w_m \bm y_j + b_m)) + \alpha_0  - f(\bm y_j)\right]^2 \\
&\;\;\quad + \frac{1}{J_2}\sum_{j=1}^{J_2}\left[\sum_{m=1}^M\alpha_m\, \mathcal{B}(\sigma( \bm w_m \bm y_j + b_m)) + \alpha_0  - g(\bm y_j)\right]^2\Bigg\}
 \end{aligned}
\end{equation}
where the parameters $\{\bm w_m\}_{m=1}^M$ and $\{b_m\}_{m=1}^M$ are first computed using the strategy in Section \ref{sec:uniform} and Eq.~\eqref{eq:rep}.

{\bf Nonlinear PDEs.}\; 
When one or both operators, $\mathcal{L}$ and $\mathcal{B}$, are nonlinear, there are two approaches to handle the situation. The first way is to wrap the least squares problem with a well established nonlinear iterative solver, e.g., Picard's methods, to solve the PDE. Within each iteration, the PDE is linearized such that we can update the coefficient $\bm \alpha$ by solving the least squares problem as mentioned above. When there is sufficient knowledge to choose a proper nonlinear solver, we prefer this approach because the well-established theory on nonlinear solvers can ensure a good convergence rate. Thus, we in fcat adopt this approach for numerical experiments in this paper. 
The second feasible approach is to wrap a gradient descent optimizer around the total loss
$L(\bm y) = \|\mathcal{L}(u(\bm y)) - \mathcal{L}(u_{\rm NN}(\bm y))\|_2^2 + \|\mathcal{B}(u(\bm y)) - \mathcal{B}(u_{\rm NN}(\bm y)) \|_2^2$. Because the neural feature space $\mathcal{P}_{\rm NN}$ is fixed, the optimization will be simpler than training the entire neural network from scratch. This approach is easier to implement and suitable for scenarios that standard nonlinear solvers do not provide a satisfactory solution.
\vspace{-0.5cm}
\begin{rem}[Not using PDE's solution data]
In this work, we do not rely on any measurement data of the solution $u(\bm y)$ when using TransNet to solve PDEs, because the operators $\mathcal{L}$ and $\mathcal{B}$ in Eq.~\eqref{eq:pde} are sufficient to ensure the existence and uniqueness of the PDE's solution. On the other hand, if any extra data of $u(\bm y)$ are available, TransNet can easily incorporate it into the least squares problem in Eq.~\eqref{eq:ls} as a supervised learning loss.
\end{rem}
\subsection{Complexity and accuracy of TransNet}\label{sec:accuracy}
The complexity of TransNet is greatly reduced compared to the scenario of using SGD to train the entire network. The construction of the neural feature space $\mathcal{P}_{\rm NN}$ only involves random number generations and a simple one-dimensional optimization in Eq.~\eqref{eq:ls2}. Moreover, these cost are completely offline, and the constructed $\mathcal{P}_{\rm NN}$ is transferable to various PDE problems. The online operation for solving linear PDEs only 
requires solving one least squares problem, where the assembling of the least squares matrix can be efficiently done using the autograd function in Tensorflow or Pytorch. The numerical experiments in Section \ref{sec:exp} show that that the accuracy and efficiency of TransNet is significantly improved compared with several baseline methods, because our method does not suffer from the 
slow convergence of SGD in neural network training.

\section{Numerical experiments}\label{sec:exp}
%
We now demonstrate the performance of TransNet by testing several classic steady-state or time-dependent PDEs in two and three dimensional spaces. In Section \ref{sec:basis}, we illustrate how to construct the transferable feature space $\mathcal{P}_{\rm NN}$. To test and demonstrate the transferability of our model, we build and test two neural features spaces, one for the 2D case and the other for the 3D case\footnote{Note that the dimension of the feature space is the sum of both space and time dimensions since it doesn't differ them.}.
The constructed feature spaces are then used  in Section \ref{sec:pde2} to solve the model PDE problems. 

\subsection{Uniform neuron distribution}
\label{sec:basis}
This experiment is to use and test the algorithm proposed in Section \ref{sec:par} to construct  transferable neural feature spaces $\mathcal{P}_{\rm NN}$ in the 2D and 3D unit balls. 
We tune the shape parameter 
$\gamma = \gamma_m$ for $m = 1, \ldots, M$ in Eq.~\eqref{eq:plane2} with $K =50$ realizations of the Gaussian process.  In addition, we also test the effect of the correlation length and the number of hidden neurons by setting different values for $\eta$ and $M$.  For each setting of $\eta$ and $M$, the shape parameter $\gamma$ is tuned separately. Additional information about the experiment setup is given in Appendix \ref{app:gp}.

\vspace{-0.2cm}
\begin{figure}[h!]
  \centering  \includegraphics[width=0.45\textwidth]{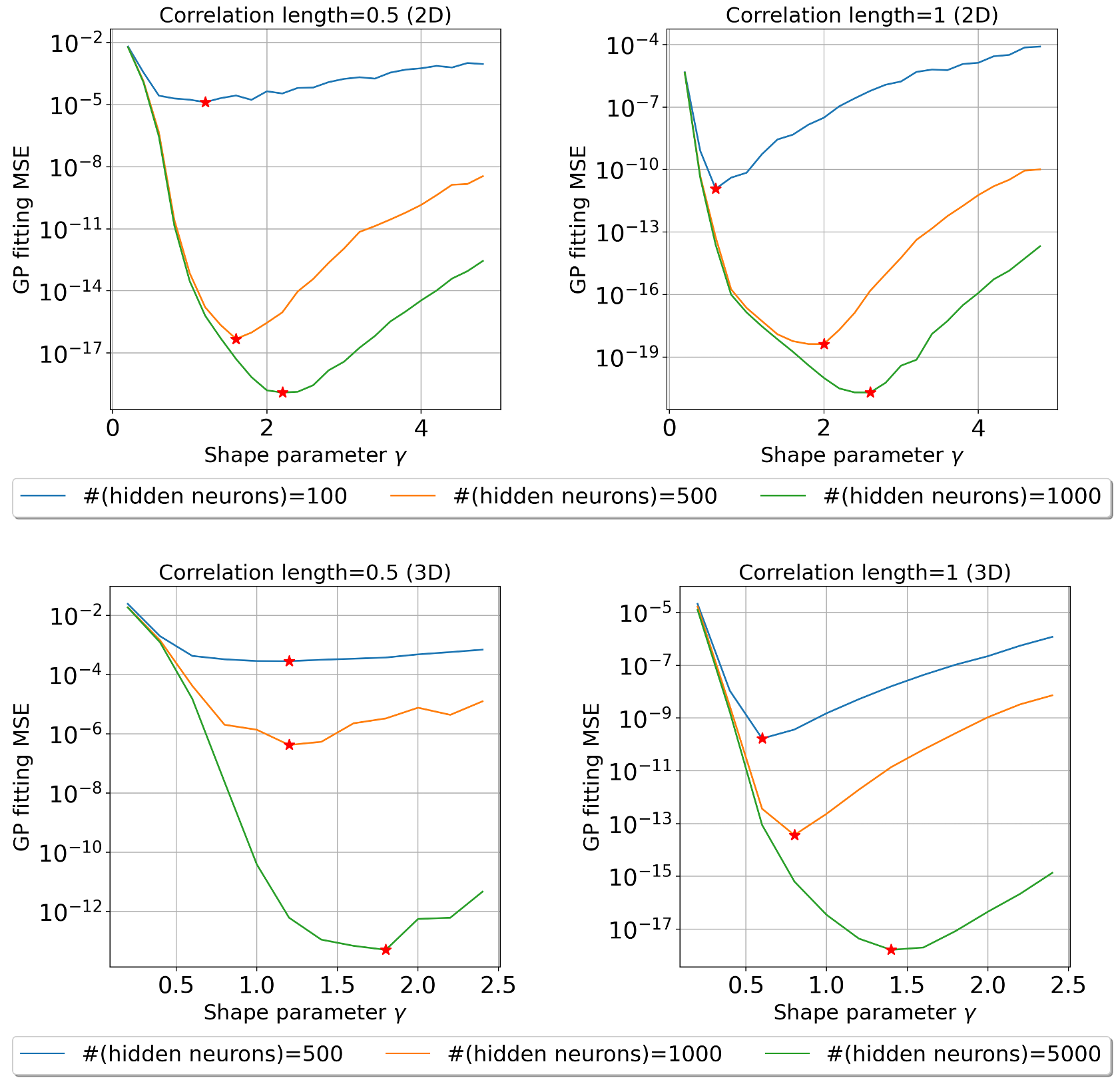}
  \vspace{-0.0cm}
      \caption{The loss landscapes of the optimizing problem in Eq.~\eqref{eq:ls2} for tuning the shape parameter $\gamma$ of the feature space $\mathcal{P}_{\rm NN}$ in two and three dimensional cases. The blue star is the optimal value for $\gamma$ founded by our method. It shows that the optimal value for $\gamma$ varies with the number of hidden neurons, meaning that tuning $\gamma$ is a necessary operation to achieve  optimal accuracy of $u_{\rm NN}$ when changing the number of hidden neurons. 
      }\label{fig:gamma}
\end{figure}

Figure \ref{fig:gamma} illustrates the landscapes of the loss function $\sum_{k=1}^K \min_{\bm \alpha}[\text{MSE}({u}_{\rm NN}(\bm y), G(\bm y| \omega_k, \eta)) ]$ of the optimization problem in Eq.~\eqref{eq:ls2} for 2D and 3D neural feature spaces. We report the results for two correlation lengths ($\eta = 0.5$ and $\eta = 1.0$) combined with three numbers of hidden neurons ($M = 100, 500, 1000$ for 2D  and $M = 500, 1000, 5000$ for 3D). We observe that the loss function behaves roughly like a parabolic curve for a fixed number of hidden neurons, so that the problem in Eq.~\eqref{eq:ls2} can be solved by a simple solver for one-dimensional optimization. More importantly, we observe that the optimal value for $\gamma$ varies with the number of hidden neurons. This provides an important insight that tuning $\gamma$ is a necessary operation to achieve optimal accuracy of $u_{\rm NN}$ when changing the number of hidden neurons. 

\begin{figure}[h!]
  \centering
  {\includegraphics[width=0.48\textwidth]{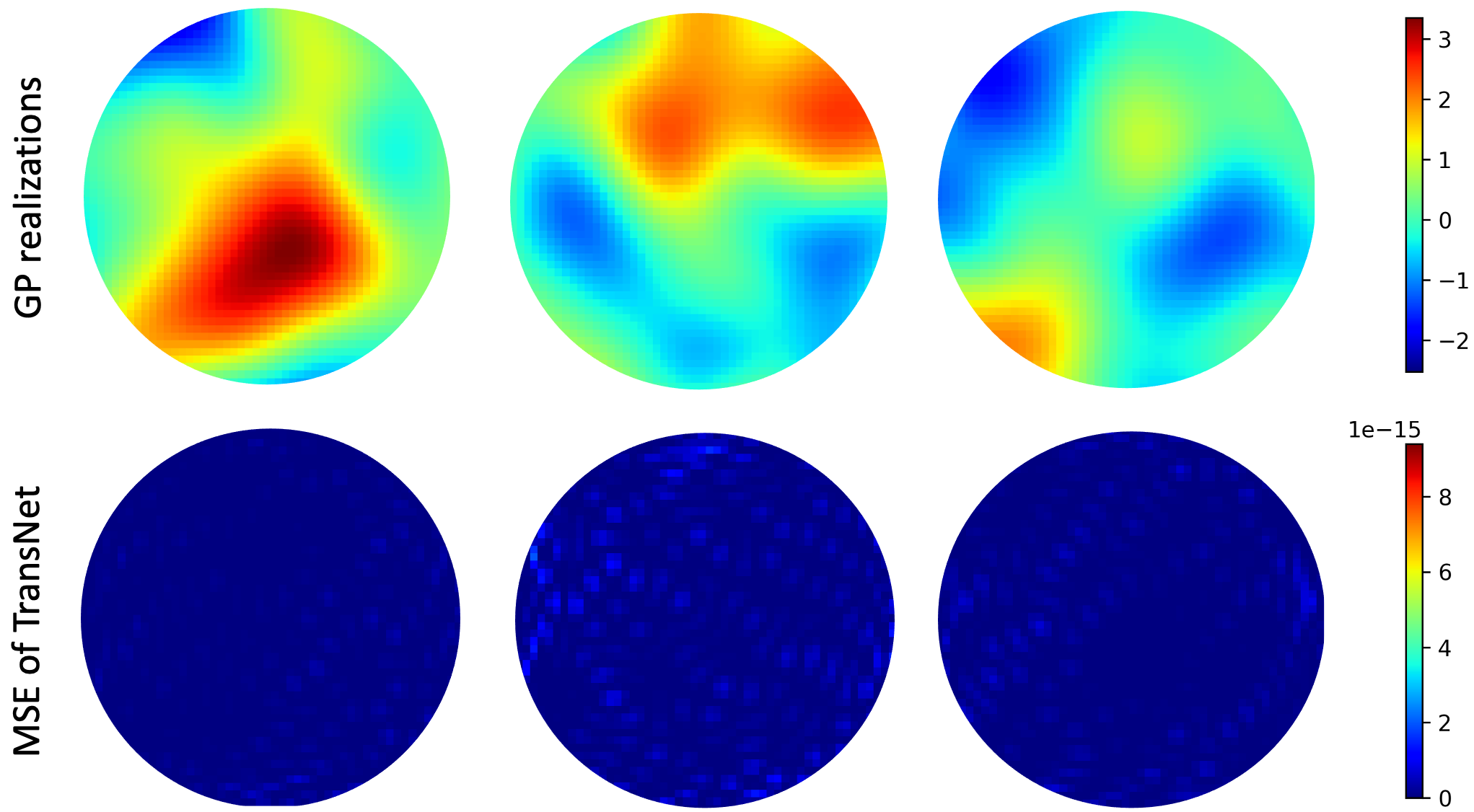}}
  \vspace{-0.2cm}
  \caption{Top row: three realizations of the auxiliary Gaussian process with the correlation length $\eta = 0.5$. Bottom row: the distribution of the MSE of TransNet's approximation with 1000 hidden neurons. Thanks to the feature space with the uniform density in the 2D unit ball (illustrated in Figure {\color{blue}\ref{fig1}(c)}), we obtain a TransNet approximation with very small MSE fluctuation.}\label{fig:gp-err}
  \vspace{-0.4cm}
\end{figure}

\begin{figure*}[h!]
    \centering
    \includegraphics[width=0.98\textwidth]{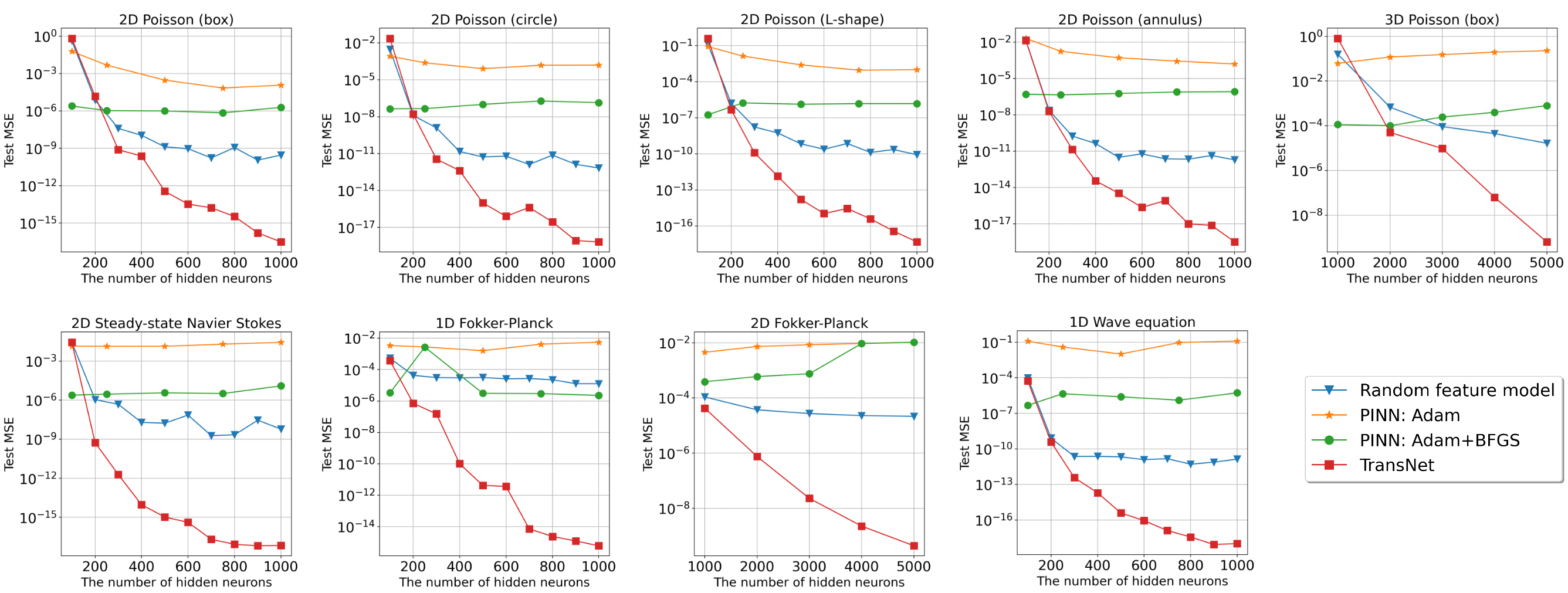}
    \caption{The MSE decay along with the increasing of the number of hidden neurons for $(C_1)$ to $(C_9)$, where all the methods use the same network architecture. Our TransNet significantly outperforms the baseline methods from two aspects: (i) {\em Transferability}: for a fixed number of hidden neurons, TransNet only need use one 2D feature space and one 3D feature space; (ii) {\em Accuracy}: TransNet achieves several orders of magnitude smaller MSE than PINN and the random feature models. TransNet does not suffer from the slow convergence in SGD-based neural network training, and can exploit more expressive power of a given neural network $u_{\rm NN}$ to obtain more accurate PDE solutions.}
    \label{fig:err}
    \vspace{-0.0cm}
\end{figure*}

Figure \ref{fig:gp-err} illustrates the error distribution when using TransNet  to approximate three realizations of the Gaussian process with correlation length $\eta = 0.5$ in the 2D unit ball. Even though the purpose of TransNet is not to approximate the Gaussian process, it is interesting to check whether the uniform density $D_M(\bm y)$ (proved in Theorem \ref{thm1}) leads to uniform error distribution.  
We use 1000 hidden neurons and the shape parameter $\gamma$ is set to 2. The bottom row of Figure \ref{fig:gp-err} shows that the MSE error distributes uniformly in the unit ball, which demonstrates the effectiveness of the feature space generation method proposed in Section \ref{sec:par}.

\subsection{PDE examples}\label{sec:pde2}
We then use the constructed 2D and 3D neural feature spaces from Section \ref{sec:basis} to solve two steady-state PDEs (i.e., the Poisson equation and the time-independent Navior-Stokes equation) and  two time-dependent PDEs (i.e., the Fokker-Planck equation and the wave equation). The definitions of the PDEs under consideration are given in Appendix \ref{app:pde}. We perform the following testing cases: 
\vspace{-0.2cm}
\begin{itemize}[leftmargin=30pt]\itemsep-0.1cm
    \item[($C_1$)] Poisson equation (2D space) in a box domain;
    \item[($C_2$)] Poisson equation (2D space) in a circular domain;
    \item[($C_3$)] Poisson equation (2D space) in an L-shaped domain;
    \item[($C_4$)] Poisson equation (2D space) in an annulus domain;
    \item[($C_5$)] Poisson equation (3D space) in a box domain;
    \item[($C_6$)] Steady-state Navier-Stokes equation (2D space);
    \item[($C_7$)] Fokker-Planck equation (1D space + 1D time);
    \item[($C_8$)] 2D Fokker-Planck equation (2D space + 1D time);
    \item[($C_9$)] 1D wave equation (1D space + 1D time)
\end{itemize}
\vspace{-0.3cm}
to demonstrate the transferability of TransNet in solving various PDEs in different domains. Recall that for time-dependent PDEs,  the temporal variable is simply treated as an extra dimension, so that we will use the 2D feature space to solve problems $(C_7)$ and $(C_9)$ and  the 3D feature space to solve problem $(C_8)$. 
We compare our method with two baseline methods, i.e., the random feature mode and the PINN. All the methods use the same network architecture, i.e., Eq.~\eqref{eq:fc} with the $tanh$ activation. Additional information about the setup of the experiments are given in Appendix \ref{app:setup}.
\begin{table*}[h!]
\footnotesize
    \centering
    \begin{tabular}{c|c|c|c|c|c|c|c|c|c}
    \toprule
     & $(C_1)$ & $(C_2)$ & $(C_3)$ & $(C_4)$ & $(C_5)$ & $(C_6)$ & $(C_7)$ & $(C_8)$ & $(C_9)$ \\
     \midrule
     Random feature model & 0.25s & 0.22s & 0.22s & 0.19s & 0.96s & 12.85s & 0.92s & 1.21s & 0.47s  \\
     PINN:Adam & 29.69s & 25.34s & 24.57s & 22.24s & 110.59s & 69.73s & 61.45s & 97.12s & 49.25s\\
     PINN:Adam+BFGS & 125.78s & 121.46s&  120.93s & 119.24s &  264.62s & 191.53s & 172.86s & 178.99s & 152.71s\\
     TransNet & 0.27s & 0.20s & 0.20s & 0.17s & 1.03s & 11.14s & 0.97s & 1.27s & 0.51s  \\
     \toprule
    \end{tabular}
    \vspace{-0.2cm}
    \caption{The computing times of TransNet and the baselines in solving the nine PDE test cases with 1000 hidden neurons. TransNet and the random feature model are significantly faster than PINN because SGD is not required in them.}
    \label{tab1}
    \vspace{-0.2cm}
\end{table*}

Figure \ref{fig:err} shows the MSE decay with the increasing of the number of the hidden neurons, where the number of hidden neurons are chosen as $M =$ 100, 200, 300, 400, 500, 600, 700, 800, 900, 1000, respectively,  for the 2D feature space, and $M =$ 1000, 2000, 3000, 4000, 5000, respectively, for the 3D feature space. 
We observe that our TransNet achieves a superior performance for all the nine test cases, which demonstrates the outstanding  transferability of TransNet.
PINN with BFGS acceleration provides a good accuracy gain compared with PINN with Adam, which means the landscape of the PDE loss exhibits severe ill-conditioning as the SGD method approaches the minimizer\footnote{BFGS can alleviate ill-conditioning by exploiting the second-order information, e.g., the approximate Hessian.}. In comparison, TransNet does not require SGD in solving the PDEs, so that TransNet does not suffer from the slow convergence of SGD used in PINN.

\begin{figure}[h!]
    \centering \includegraphics[width=0.45\textwidth]{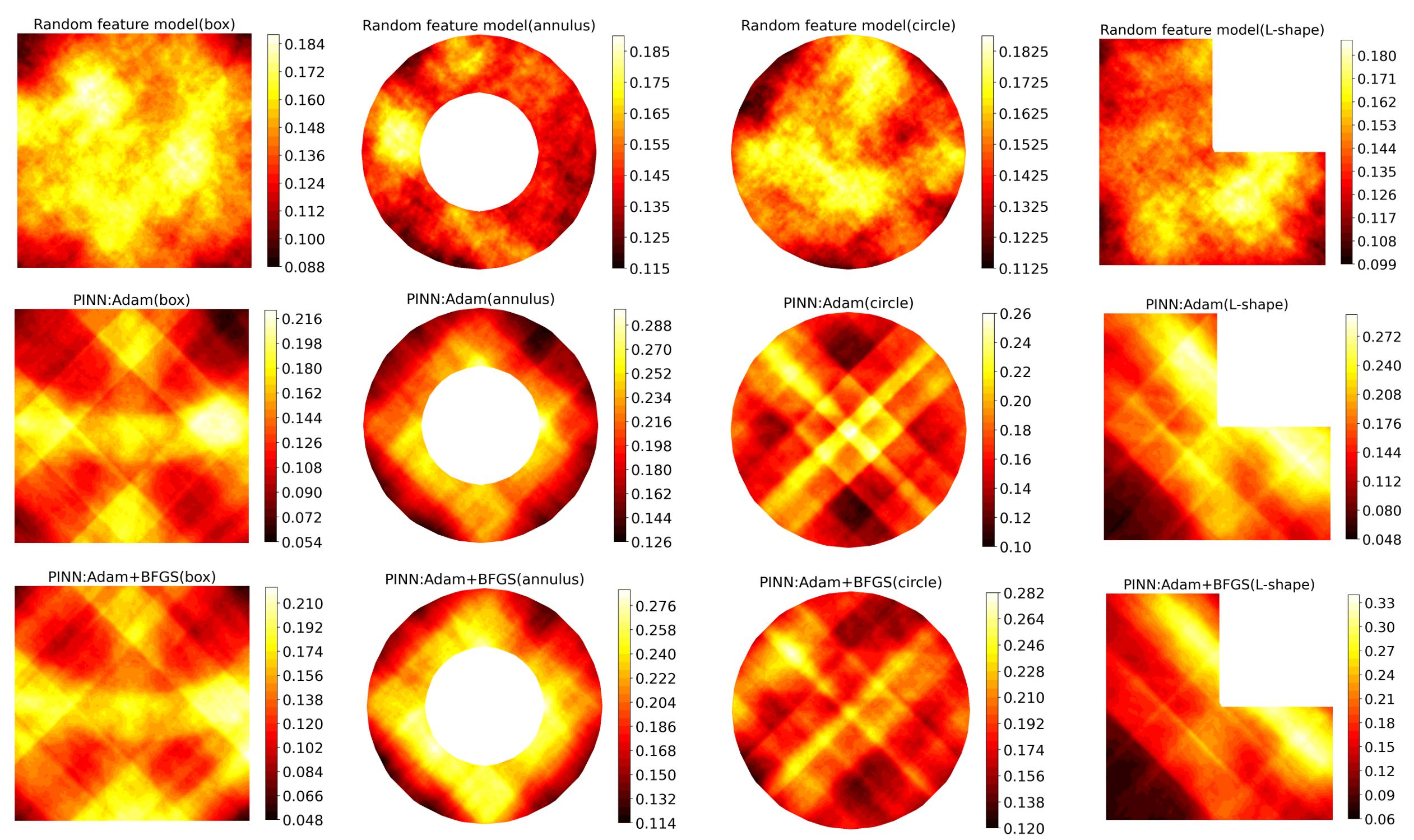}
    \vspace{-0.2cm}
    \caption{The density function $D_M(\bm y)$ with $\tau =0.2$ in Eq.~\eqref{eq:den} of the neural feature spaces obtained by training PINN and the random feature models in solving the Poisson equation in the 2D space, i.e., problems $(C_1)$ - $(C_4)$. Compared to the uniform density of TransNet in Figure \ref{fig1}{\color{blue}(c)}, both PINN and the random feature model cannot provide  feature spaces with  uniform density, which is one explanation of their under-performance shown in Figure \ref{fig:err}.}
    \label{fig:PINN_den}
    \vspace{-0.1cm}
\end{figure}

Figure \ref{fig:PINN_den} shows the density function $D_M(\bm y)$ in Eq.~\eqref{eq:den} of the feature spaces obtained by training PINN and the random feature models in solving the Poisson equation in the 2D space, i.e., case $(C_1)$ - $(C_4)$, where the constant $\tau$ in Eq.~\eqref{eq:den} is set to 0.2. Compared with TransNet's uniform density shown in Figure {\ref{fig1}\color{blue}(c)}, the feature spaces obtained by the baseline methods have highly non-uniform densities in the domain of computation. The random feature models tend to have higher density, i.e., more hidden neurons, near the center of the domain. 
The first row in Figure \ref{fig:PINN_den} can be viewed as the initial densities of the feature space for PINN; the second and the third rows are the final densities. We can see that the training of PINN does not necessarily lead to a more uniform density function $D_M(\bm y)$, which is one of the reasons why PINN cannot exploit the full expressive power of the neural network $u_{\rm NN}$. 
\vspace{-0.5cm}

\vspace{0.5cm}
\section{Conclusion}\label{sec:con}
We propose a transferable neural network model to advance the state of the art of using neural networks to solve PDEs. The key ingredient is to construct a neural feature space independent of any PDE, which makes it easy to transfer the neural feature space to various PDEs in different domains. Moreover, because the feature space is in fact fixed when using TransNet to solve a PDE, we only need to solve linear least squares problems, which avoids the drawbacks of SGD-based training algorithms, e.g., ill-conditioning. Numerical experiments show that the proposed TransNet can exploit more expressive power of a given neural network than the compared baselines. This work is the first scratch in this research direction, and there are multiple potential related topics that will be studied in our future work, including (1) {\em theoretical analysis of the convergence rate of TransNet in solving PDEs.} We observe in Figure \ref{fig:err} that the MSE of TransNet decays along with the increasing of the number of hidden neurons. A natural question to study is that whether TransNet  can achieve the optimal convergence rate of the 
single-hidden-layer fully-connected neural network. (2) {\em Extension to multi-layer neural networks.} Even though the single-hidden-layer model has sufficient expressive power for the PDEs tested in this work, there are more complicated PDEs, e.g., turbulence models, that could require multi-layer models with much higher expressive power. (3) {\em The properties of the least squares problem.} In this work, we use the standard least squares solver of Pytorch in the numerical experiments. However, it is worth further investigation of the properties of this specific least squares problem. For example, since the set of neurons $\{\sigma(\bm w_m \bm y + b_m)\}_{m=1}^M$ forms a non-orthogonal basis, it is possible to have linearly correlated neurons which will reduce the column rank of the least squares matrix, or even lead to an under-determined system. This will require the use of some regularization techniques, e.g., ridge regression, to stabilize the least squares system. Additionally, compressed sensing, i.e., $\ell_1$ regularization, could be added  to remove redundant neurons from the feature space as needed and obtain a sparse neural network.

\section*{Acknowledgement}
This work was supported by the U.S. Department of Energy, Office of Science, Office of Advanced Scientific Computing Research, Applied Mathematics Program, under the contract number ERKJ387. This work was accomplished at Oak Ridge National Laboratory (ORNL). ORNL is operated by UT-Battelle, LLC., for the U.S. Department of Energy under Contract DE-AC05-00OR22725. 


\begin{thebibliography}{33}
\providecommand{\natexlab}[1]{#1}
\providecommand{\url}[1]{\texttt{#1}}
\expandafter\ifx\csname urlstyle\endcsname\relax
  \providecommand{\doi}[1]{doi: #1}\else
  \providecommand{\doi}{doi: \begingroup \urlstyle{rm}\Url}\fi

\bibitem[Anitescu et~al.(2019)Anitescu, Atroshchenko, Alajlan, and
  Rabczuk]{anitescu2019artificial}
Anitescu, C., Atroshchenko, E., Alajlan, N., and Rabczuk, T.
\newblock Artificial neural network methods for the solution of second order
  boundary value problems.
\newblock \emph{Computers, Materials \& Continua}, 59\penalty0 (1):\penalty0
  345--359, 2019.

\bibitem[Bach(2017)]{10.5555/3122009.3122030}
Bach, F.
\newblock On the equivalence between kernel quadrature rules and random feature
  expansions.
\newblock \emph{J. Mach. Learn. Res.}, 18\penalty0 (1):\penalty0 714–751, jan
  2017.
\newblock ISSN 1532-4435.

\bibitem[Chakraborty(2020)]{Chakraborty2020TransferLB}
Chakraborty, S.~L.
\newblock Transfer learning based multi-fidelity physics informed deep neural
  network.
\newblock \emph{J. Comput. Phys.}, 426:\penalty0 109942, 2020.

\bibitem[Chen \& Chen(1995)Chen and Chen]{chen1995universal}
Chen, T. and Chen, H.
\newblock Universal approximation to nonlinear operators by neural networks
  with arbitrary activation functions and its application to dynamical systems.
\newblock \emph{IEEE Transactions on Neural Networks}, 6\penalty0 (4):\penalty0
  911--917, 1995.

\bibitem[{Clark Di Leoni} et~al.(2023){Clark Di Leoni}, Lu, Meneveau,
  Karniadakis, and Zaki]{CLARKDILEONI2023111793}
{Clark Di Leoni}, P., Lu, L., Meneveau, C., Karniadakis, G.~E., and Zaki, T.~A.
\newblock Neural operator prediction of linear instability waves in high-speed
  boundary layers.
\newblock \emph{Journal of Computational Physics}, 474:\penalty0 111793, 2023.

\bibitem[Desai et~al.(2021)Desai, Mattheakis, Joy, Protopapas, and
  Roberts]{Desai2021OneShotTL}
Desai, S., Mattheakis, M., Joy, H., Protopapas, P., and Roberts, S.~J.
\newblock One-shot transfer learning of physics-informed neural networks.
\newblock \emph{ArXiv}, abs/2110.11286, 2021.

\bibitem[Dissanayake \& Phan-Thien(1994)Dissanayake and
  Phan-Thien]{dissanayake1994neural}
Dissanayake, M. and Phan-Thien, N.
\newblock Neural-network-based approximations for solving partial differential
  equations.
\newblock \emph{Communications in Numerical Methods in Engineering},
  10\penalty0 (3):\penalty0 195--201, 1994.

\bibitem[E \& Yu(2018)E and Yu]{e2018ritz}
E, W. and Yu, B.
\newblock The deep {Ritz} method: a deep learning-based numerical algorithm for
  solving variational problems.
\newblock \emph{Communication in Mathematics and Statistics}, 6\penalty0
  (1):\penalty0 1--12, 2018.

\bibitem[Gin et~al.(2021)Gin, Shea, Brunton, and Kutz]{gin2021deepgreen}
Gin, C.~R., Shea, D.~E., Brunton, S.~L., and Kutz, J.~N.
\newblock Deepgreen: deep learning of green’s functions for nonlinear
  boundary value problems.
\newblock \emph{Scientific Reports}, 11\penalty0 (1):\penalty0 1--14, 2021.

\bibitem[Jin et~al.(2022)Jin, Meng, and Lu]{JML2022MIONet}
Jin, P., Meng, S., and Lu, L.
\newblock Mionet: Learning multiple-input operators via tensor product.
\newblock \emph{SIAM Journal on Scientific Computing}, 44\penalty0
  (6):\penalty0 A3490--A3514, 2022.

\bibitem[{Karniadakis} et~al.(2021){Karniadakis}, {Kevrekidis}, {Lu},
  {Perdikaris}, {Wang}, and {Yang}]{2021NatRP...3..422K}
{Karniadakis}, G.~E., {Kevrekidis}, I.~G., {Lu}, L., {Perdikaris}, P., {Wang},
  S., and {Yang}, L.
\newblock {Physics-informed machine learning}.
\newblock \emph{Nature Reviews Physics}, 3\penalty0 (6):\penalty0 422--440,
  January 2021.
\newblock \doi{10.1038/s42254-021-00314-5}.

\bibitem[Krishnapriyan et~al.(2021)Krishnapriyan, Gholami, Zhe, Kirby, and
  Mahoney]{krishnapriyan2021characterizing}
Krishnapriyan, A., Gholami, A., Zhe, S., Kirby, R., and Mahoney, M.~W.
\newblock Characterizing possible failure modes in physics-informed neural
  networks.
\newblock \emph{Advances in Neural Information Processing Systems}, 34, 2021.

\bibitem[Lagaris et~al.(1998)Lagaris, Likas, and
  Fotiadis]{lagaris1998artificial}
Lagaris, I.~E., Likas, A., and Fotiadis, D.~I.
\newblock Artificial neural networks for solving ordinary and partial
  differential equations.
\newblock \emph{IEEE transactions on neural networks}, 9\penalty0 (5):\penalty0
  987--1000, 1998.

\bibitem[Li et~al.(2020)Li, Kovachki, Azizzadenesheli, Liu, Stuart,
  Bhattacharya, and Anandkumar]{li2020multipole}
Li, Z., Kovachki, N., Azizzadenesheli, K., Liu, B., Stuart, A., Bhattacharya,
  K., and Anandkumar, A.
\newblock Multipole graph neural operator for parametric partial differential
  equations.
\newblock \emph{Advances in Neural Information Processing Systems},
  33:\penalty0 6755--6766, 2020.

\bibitem[Li et~al.(2021{\natexlab{a}})Li, Kovachki, Azizzadenesheli,
  Bhattacharya, Stuart, Anandkumar, et~al.]{li2020fourier}
Li, Z., Kovachki, N.~B., Azizzadenesheli, K., Bhattacharya, K., Stuart, A.,
  Anandkumar, A., et~al.
\newblock Fourier neural operator for parametric partial differential
  equations.
\newblock In \emph{International Conference on Learning Representations},
  2021{\natexlab{a}}.

\bibitem[Li et~al.(2021{\natexlab{b}})Li, Zheng, Kovachki, Jin, Chen, Liu,
  Azizzadenesheli, and Anandkumar]{li2021physics}
Li, Z., Zheng, H., Kovachki, N., Jin, D., Chen, H., Liu, B., Azizzadenesheli,
  K., and Anandkumar, A.
\newblock Physics-informed neural operator for learning partial differential
  equations.
\newblock \emph{arXiv preprint arXiv:2111.03794}, 2021{\natexlab{b}}.

\bibitem[Liu et~al.(2022{\natexlab{a}})Liu, Huang, Chen, and Suykens]{9495136}
Liu, F., Huang, X., Chen, Y., and Suykens, J. A.~K.
\newblock Random features for kernel approximation: A survey on algorithms,
  theory, and beyond.
\newblock \emph{IEEE Transactions on Pattern Analysis and Machine
  Intelligence}, 44\penalty0 (10):\penalty0 7128--7148, 2022{\natexlab{a}}.

\bibitem[Liu et~al.(2022{\natexlab{b}})Liu, McCalla, and
  Schaeffer]{https://doi.org/10.48550/arxiv.2212.05591}
Liu, Y., McCalla, S.~G., and Schaeffer, H.
\newblock Random feature models for learning interacting dynamical systems,
  2022{\natexlab{b}}.

\bibitem[Long et~al.(2018)Long, Lu, Ma, and Dong]{long2018pde}
Long, Z., Lu, Y., Ma, X., and Dong, B.
\newblock {PDE-Net: Learning PDEs from data}.
\newblock In \emph{International Conference on Machine Learning}, pp.\
  3214--3222, 2018.

\bibitem[Long et~al.(2019)Long, Lu, and Dong]{long2019pde}
Long, Z., Lu, Y., and Dong, B.
\newblock {PDE-Net} 2.0: Learning {PDEs} from data with a numeric-symbolic
  hybrid deep network.
\newblock \emph{Journal of Computational Physics}, 399:\penalty0 108925, 2019.

\bibitem[Lu et~al.(2021{\natexlab{a}})Lu, Jin, Pang, Zhang, and
  Karniadakis]{lu2021learning}
Lu, L., Jin, P., Pang, G., Zhang, Z., and Karniadakis, G.~E.
\newblock Learning nonlinear operators via {DeepONet} based on the universal
  approximation theorem of operators.
\newblock \emph{Nature Machine Intelligence}, 3\penalty0 (3):\penalty0
  218--229, 2021{\natexlab{a}}.

\bibitem[Lu et~al.(2021{\natexlab{b}})Lu, Meng, Mao, and
  Karniadakis]{lu2021deepxde}
Lu, L., Meng, X., Mao, Z., and Karniadakis, G.~E.
\newblock Deepxde: A deep learning library for solving differential equations.
\newblock \emph{SIAM Review}, 63\penalty0 (1):\penalty0 208--228,
  2021{\natexlab{b}}.

\bibitem[McDonald \& \'{A}lvarez(2021)McDonald and
  \'{A}lvarez]{NEURIPS2021_72fe6f9f}
McDonald, T. and \'{A}lvarez, M.
\newblock Compositional modeling of nonlinear dynamical systems with ode-based
  random features.
\newblock In Ranzato, M., Beygelzimer, A., Dauphin, Y., Liang, P., and Vaughan,
  J.~W. (eds.), \emph{Advances in Neural Information Processing Systems},
  volume~34, pp.\  13809--13819. Curran Associates, Inc., 2021.

\bibitem[Nelsen \& Stuart(2021)Nelsen and Stuart]{doi:10.1137/20M133957X}
Nelsen, N.~H. and Stuart, A.~M.
\newblock The random feature model for input-output maps between banach spaces.
\newblock \emph{SIAM Journal on Scientific Computing}, 43\penalty0
  (5):\penalty0 A3212--A3243, 2021.

\bibitem[Quarteroni et~al.(2007)Quarteroni, Sacco, and
  Saleri]{2013JSV...332.4403B}
Quarteroni, A., Sacco, R., and Saleri, F.
\newblock \emph{{Numerical Mathematics}}, volume 332.
\newblock Springer Science Business Media {\&}, 2007.

\bibitem[Raissi et~al.(2019)Raissi, Perdikaris, and
  Karniadakis]{raissi2019physics}
Raissi, M., Perdikaris, P., and Karniadakis, G.~E.
\newblock Physics-informed neural networks: A deep learning framework for
  solving forward and inverse problems involving nonlinear partial differential
  equations.
\newblock \emph{Journal of Computational Physics}, 378:\penalty0 686--707,
  2019.

\bibitem[Sirignano \& Spiliopoulos(2018)Sirignano and
  Spiliopoulos]{sirignano2018}
Sirignano, J. and Spiliopoulos, K.
\newblock {DGM}: A deep learning algorithm for solving partial differential
  equations.
\newblock \emph{Journal of Computational Physics}, 375:\penalty0 1339--1354,
  2018.

\bibitem[Sun et~al.(2018)Sun, Gilbert, and Tewari]{Sun2018OnTA}
Sun, Y., Gilbert, A.~C., and Tewari, A.
\newblock On the approximation capabilities of relu neural networks and random
  relu features.
\newblock \emph{arXiv: Machine Learning}, 2018.

\bibitem[Teng et~al.(2022)Teng, Zhang, Wang, and Ju]{Teng2022gfnet}
Teng, Y., Zhang, X., Wang, Z., and Ju, L.
\newblock Learning green's functions of linear reaction-diffusion equations
  with application to fast numerical solver.
\newblock In \emph{Mathematical and Scientific Machine Learning Conference},
  2022.

\bibitem[Wang et~al.(2021)Wang, Wang, and Perdikaris]{wang2021learning}
Wang, S., Wang, H., and Perdikaris, P.
\newblock Learning the solution operator of parametric partial differential
  equations with physics-informed deeponets.
\newblock \emph{Science Advances}, 7\penalty0 (40):\penalty0 eabi8605, 2021.

\bibitem[Zang et~al.(2020)Zang, Bao, Ye, and Zhou]{zang2019weak}
Zang, Y., Bao, G., Ye, X., and Zhou, H.
\newblock Weak adversarial networks for high dimensional partial differential
  equations.
\newblock \emph{Journal of Computational Physics}, 411:\penalty0 109409, 2020.

\bibitem[Zhang et~al.(2021)Zhang, Cheng, and Ju]{zhang2021ifnn}
Zhang, X., Cheng, T., and Ju, L.
\newblock Implicit form neural network for learning scalar hyperbolic
  conservation laws.
\newblock In \emph{Mathematical and Scientific Machine Learning Conference},
  pp.\  1082--1098, 2021.

\bibitem[Zhao \& Wright(2021)Zhao and Wright]{zhao2020solving}
Zhao, J. and Wright, C.~L.
\newblock Solving allen-cahn and cahn-hilliard equations using the adaptive
  physics informed neural networks.
\newblock \emph{Communications in Computational Physics}, 29:\penalty0
  930--954, 2021.

\end{thebibliography}

\newpage
\onecolumn
\section*{Appendix}
\appendix

\section{The proof of Theorem \ref{thm1}}\label{sec:proof}
For the re-parameterization in Eq.~\eqref{eq:plane2}, we can treat $\{\bm a_m\}_{m=1}^M$ as $M$ independent and identically distributed (i.i.d.) random variables on the $d$-dimensional unit sphere, and $\{r_m\}_{m=1}^M$
as $M$ i.i.d. random variables following the uniform distribution in $[0,1]$. For a fixed $\bm y \in \Omega$, the expectation of $D_M(\bm y)$ is 
\begin{equation}\label{eq:bb}
\mathbb{E}[D_M(\bm y)] =  \frac{1}{M}\sum_{m=1}^M \mathbb{E}\left[\mathbf{1}_{dist(\bm y, m) < \tau} (\bm y)\right].
\end{equation}
Because $\mathbb{E}[\mathbf{1}_{dist(\bm y, m) < \tau} (\bm y)] = \mathbb{E}[\mathbf{1}_{dist(\bm y, m') < \tau} (\bm y)]$, we only need to calculate one expectation $\mathbb{E}[\mathbf{1}_{dist(\bm y, m) < \tau} (\bm y)]$. Therefore, we can drop the subscript of $\bm a_m$ and use $\bm a$ to denote $\bm a_m$ in the following derivation. 

To proceed, we define the representations of the vectors $\bm a = (a_{1}, \ldots, a_{d})$ and $\bm y = (y_1, \ldots, y_d)$ under different coordinate systems. We denote by $\mathcal{C}_{\rm original}$ the original Cartesian coordinate system and denote by $\bm a | \mathcal{C}_{\rm original}$ and $\bm y|\mathcal{C}_{\rm original}$ the representations of $\bm a$ and $\bm y$ under $\mathcal{C}_{\rm original}$. Because $\bm a $ and $\bm y$ are defined in $\mathcal{C}_{\rm original}$, we have 
\[
\bm a | \mathcal{C}_{\rm original} = (a_{1}, \ldots, a_{d})\; \text{ and } \;
\bm y|\mathcal{C}_{\rm original} = (y_1, \ldots, y_d).
\]
We can also define a rotated Cartesian coordinate system, denoted by $\mathcal{C}_{\rm rot}$, such that the first coordinate axis of $\mathcal{C}_{\rm rot}$ aligns with the direction of $\bm y$. We denote by $\bm c_1, \ldots, \bm c_d$ the directions of the coordinate axes of $\mathcal{C}_{\rm rot}$, so the vector $\bm a$ can be represented in $\mathcal{C}_{\rm rot}$ as 
\[
\bm a | \mathcal{C}_{\rm rot} = (\tilde{a}_{1}, \ldots, \tilde{a}_{d})\; \text{ and } \; \bm a = \tilde{a}_{1} \bm c_1 + \cdots + \tilde{a}_{d} \bm c_d.
\]
Because $\bm c_1 = \bm y/\|\bm y\|_2$, we have 
\[
{\bm y} |\mathcal{C}_{\rm rot} = (\|\bm y\|_2, 0, \cdots, 0).
\]
Based on $\mathcal{C}_{\rm rot}$, we define a $d$-dimensional hyperspherical coordinate system, denoted by $\mathcal{S}_{\rm rot}$, with one radial variable $r$, $d-2$ polar angles  $(\phi_1, \ldots, \phi_{d-2})$ ranging over $[0, \pi]$ and one azimuthal angle $\phi_{d-1}$ ranging over $[0, 2\pi]$. Then, the unit vector $\bm a$ can be represented by the angular variables of $\mathcal{S}_{\rm rot}$, i.e., 
\[
\begin{aligned}
     \tilde{a}_{1} & = \cos(\phi_{1})\\
     \tilde{a}_{2} & = \sin(\phi_{1})\cos(\phi_{2})\\
     & \hspace{0.2cm}\vdots \\
     \tilde{a}_{d-1} & = \sin(\phi_{1})\cdots \sin(\phi_{d-2}) \cos(\phi_{d-1})\\ 
     \tilde{a}_{d} & = \sin(\phi_{1})\cdots \sin(\phi_{d-2}) \sin(\phi_{d-1}).\\ 
\end{aligned}
\]
where $(\phi_{1}, \ldots, \phi_{d-1})$ are the representation of $\bm a$ under $\mathcal{S}_{\rm rot}$. Since inner product is independent of coordinate system, the inner product $\bm a \bm y$ can be performed under $\mathcal{C}_{\rm rot}$ to obtain
\[
\bm a \bm y = (\bm a|\mathcal{C}_{\rm rot}) (\bm y |\mathcal{C}_{\rm rot})  = \|\bm y\|_2 \tilde{a}_{1} + 0\, \tilde{a}_{2} + \ldots + 0\, \tilde{a}_{d} = \|\bm y\|_2 \cos(\phi_{1}),
\]
which is independent of $\phi_{2}, \ldots, \phi_{d-1}$.

Now we derive the probability density function of the inner product $\bm a \bm y$ for a fixed $\bm y$. For any fixed $\phi_{2}, \ldots, \phi_{d-1}$, the set
\[
\mathcal{J}_{\phi_{1}|\phi_{2}, \ldots, \phi_{d-1}} := \{(1,\phi_{1}, \phi_{2}, \ldots, \phi_{d-1}) \;|\; \phi_{1} \in [0, \pi] \text{ and }  \phi_{2}, \ldots, \phi_{d-1} \text{ are fixed.}\},
\]
is a one-dimensional half circle on the $d$-dimensional unit sphere. When $\bm a$ is uniformly distributed on the $d$-dimensional unit sphere, the conditional variable $\bm a | (\phi_{2}, \ldots, \phi_{d-1})$ is uniformly distributed on the half circle $\mathcal{J}_{\phi_{1}|\phi_{2}, \ldots, \phi_{d-1}}$ and $\phi_1$ follows a uniform distribution over $[0, \pi]$ \citep{2013JSV...332.4403B}. Then, we have that the variable ${z} = \cos(\phi_1|\phi_{2}, \ldots, \phi_{d-1})$ follows the Chebyshev density
\begin{equation}\label{eq:pz}
    p_Z(z) = \frac{1}{\pi} \frac{1}{\sqrt{1-z^2}} \;\; z \in [-1,1],
\end{equation}
for any fixed $(\phi_2, \ldots, \phi_{d-1})$. 
Because the inner product $\bm a \bm y = \|\bm y\|_2 \cos(\phi_1)$ is independent of $(\phi_2, \ldots, \phi_{d-1})$, the conditional density in Eq.~\eqref{eq:pz} is also the marginal density, i.e., $p_Z(z)$ in Eq.~\eqref{eq:pz} is also the density of $z = \cos(\phi_1)$.

Next we derive the analytical form of the expectation $\mathbb{E}[\mathbf{1}_{dist(\bm y, m) < \tau} (\bm y)]$. For the convenience of derivation, we temporarily change the distribution of $r$ to a uniform distribution in $[-1,0]$, which leads to an equivalent feature space to the one with $r \in [0,1]$. 
Since $\bm a_m$ is a unit vector, we have 
$
dist(\bm y, m) = |\bm a \bm y + r|.
$
Substituting $\bm a \bm y = \|\bm y\|_2 z = \|\bm y\|_2 \cos(\phi_1)$ into $\mathbb{E}[\mathbf{1}_{dist(\bm y, m) < \tau} (\bm y)]$, we have 
\[
\begin{aligned}
\mathbb{E}[\mathbf{1}_{dist(\bm y, m) < \tau} (\bm y)] & = \mathbb{E}[\mathbf{1}_{|z\|\bm y\|_2 + r| < \tau} (\bm y)]\\
& = \int_{\{z\|\bm y\|_2 + r <\tau\} \cup \{z\|\bm y\|_2 + r >-\tau \}} p_Z(z)p_R(r) dzdr\\
& = \int_{\{z\|\bm y\|_2 + r <\tau\} \cup \{z\|\bm y\|_2 + r >-\tau \}}   \frac{1}{\pi} \frac{1}{\sqrt{1-z^2}} dzdr.\\
\end{aligned}
\]
The integral can be exactly calculated for the following two cases. 
\begin{itemize}[leftmargin=15pt]
\item {\bf Case 1}: $\|\bm y \|_2 < \tau$ meaning the integration range is below the line $r = -z\|\bm y\|_2 + \tau$. In this case, we have
\[
\begin{aligned}
\mathbb{E}[\mathbf{1}_{dist(\bm y, m) < \tau} (\bm y)] & = \mathbb{E}[\mathbf{1}_{|z\|\bm y\|_2 + r| < \tau} (\bm y)]\\
& = \int_{-1}^1 \int_{0}^{-z\|\bm y\|_2 + \tau} \frac{1}{\pi} \frac{1}{\sqrt{1-z^2}}  dr dz\\
&= \int_{-1}^1 \frac{-z\|\bm y\|_2 + \tau}{\pi \sqrt{1-z^2}}  dz\\
& =\frac{-\|\bm y\|_2}{\pi} \int_{-1}^1 \frac{ z }{ \sqrt{1-z^2}}  dz +\frac{\tau}{\pi} \int_{-1}^1 \frac{1}{ \sqrt{1-z^2}}  dz\\
& = 0 + \frac{\tau}{\pi}\left(\frac{\pi}{2} - (-\frac{\pi}{2})\right)\\
& = \tau.
\end{aligned}
\]
\item {\bf Case 2}: $\tau \le \|\bm y \|_2 \le 1- \tau$ meaning the integration range is between the line: $ r = -z\|\bm y\|_2 + \tau$ and $ r = -z\|\bm y\|_2 - \tau$. In this case, we have
\[
\begin{aligned}
\mathbb{E}[\mathbf{1}_{dist(\bm y, m) < \tau} (\bm y)] & = \mathbb{E}[\mathbf{1}_{|z\|\bm y\|_2 + r| < \tau} (\bm y)]\\
 &= \int_{-1}^0 \int_{-z\|\bm y\|_2 - \tau}^{-z\|\bm y\|_2 + \tau} \frac{1}{\pi} \frac{1}{\sqrt{1-z^2}} dr dz\\
    &= \int_{-1}^0 \frac{2\tau}{\pi \sqrt{1-z^2}}  dz\\
    &= \tau.
\end{aligned}
\]
Combining Case 1 and 2, we have 
\[
\mathbb{E}[\mathbf{1}_{dist(\bm y, m) < \tau} (\bm y)] = \tau \; \text{ for any }\; \|\bm y\|_2 \le 1-\tau.
\]
\end{itemize}
Substituting this into Eq.~\eqref{eq:bb} concludes the proof.

\section{Setup of the experiments in Section \ref{sec:basis}}\label{app:gp}
We use the python package \texttt{gstools} (https://github.com/GeoStat-Framework/GSTools/) to generate realizations of the Gaussian process.
For a fixed correlation length, we generate 10 realizations of the Gaussian process, i.e., $K = 10$ in Eq.~\eqref{eq:ls2}, to tune the shape parameter $\gamma$ of the transferable feature space. 
For the feature space for the two-dimensional PDEs, 
we sample each realization at $50^2$ uniformly distributed locations in $B_1(\bm 0)$, i.e., $J = 2500$ in Eq.~\eqref{eq:mse}, to compute the MSE in Eq.~\eqref{eq:mse}. For the feature space for the three-dimensional PDEs, we sample each realization at $50^3$, i.e., $J = 125,000$ in Eq.~\eqref{eq:mse}, to compute the MSE in Eq.~\eqref{eq:mse}. A simple grid search is used to solve the one-dimensional optimization problem in Eq.~\eqref{eq:ls2} to find the optimal shape parameter $\gamma$.

\section{Definitions of the PDEs in Section \ref{sec:pde2}}\label{app:pde}
%

The definitions of the PDEs considered in Section \ref{sec:pde2} are given below. 

\begin{figure}[h!]
    \centering
    \includegraphics[width=0.95\textwidth]{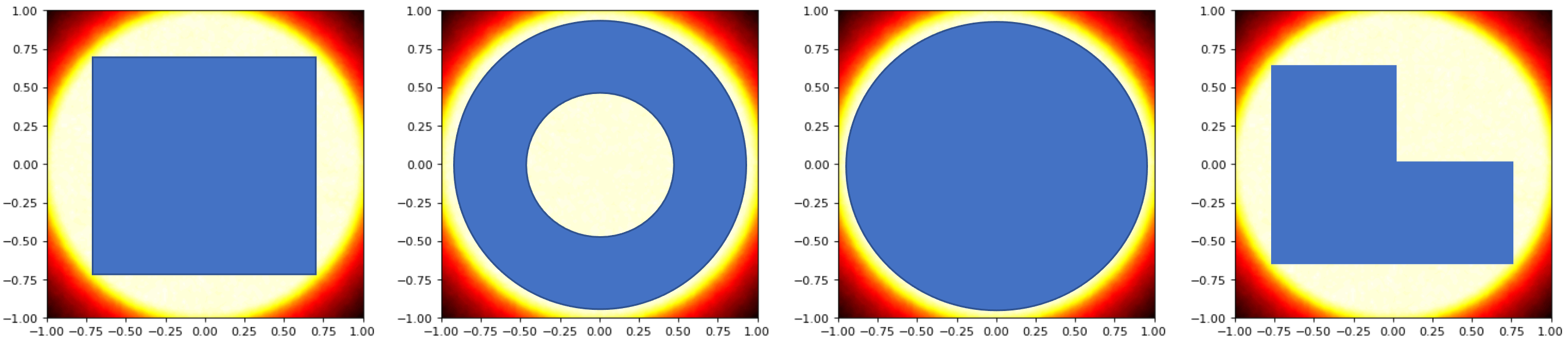}
    \vspace{-0.2cm}
    \caption{Illustration of how to place the domains of computation for the test cases $(C_1)$ -- $(C_4)$ in Section \ref{sec:pde2} into the unit ball to use the transferable feature space to solve the Poisson's equation in different domains. }
    \label{fig:shape}
\end{figure}

{\bf The Poisson's equation} considered in case $(C_1)$--$(C_5)$ is defined by
\begin{equation}
    \Delta u(\bm x) = f(\bm x),
\end{equation}
where the exact solution for the 2D settings, i.e., $(C_1)$--$(C_4)$, is 
$u(\bm x) = \sin(2\pi x_1)\sin(2\pi x_2)\sin(2\pi x_3)$, and the exact solution for the 3D setting, i.e., $(C_5)$, is $u(\bm x) = \sin(2\pi x_1)\sin(2\pi x_2)$. The forcing term $f(\bm x)$ can be obtained by applying the Laplacian operator to the exact solution. The domains of computation for $(C_1)$--$(C_5)$ are given below:
%
\begin{itemize}[leftmargin=40pt]\itemsep-0.0cm
    \item[($C_1$)] A 2D box domain: $\Omega = [-1,1]^2$; 
    \item[($C_2$)] A 2D circular domain: $\Omega = B_1(\bm 0)$;
    \item[($C_3$)] A 2D L-shaped domain: $\Omega = [-1,1]^2 \backslash [0,1]^2$;
    \item[($C_4$)] A 2D annulus domain: $\Omega = B_1(\bm 0) \backslash B_{0.5}(\bm 0)$;
    \item[($C_5$)] A 3D box domain $\Omega = [-1,1]^3$.
\end{itemize}
We consider the Dirichlet boundary condition in the experiments, where the boundary condition $g(\bm x)$ in Eq.~\eqref{eq:pde} can be obtained by restricting the exact solution on the boundary of $\Omega$. Figure \ref{fig:shape}
illustrates how to place the domains of computation into the unit ball for for the test cases $(C_1)$ -- $(C_4)$ to use the transferable feature space.

{\bf The steady-state Navier-Stokes equation} considered in case $(C_6)$ is defined by:
\[
\begin{aligned}
    {\bm u} \cdot \nabla {\bm u} + \nabla p - \nu \Delta {\bm u} & = 0\\
     \nabla \cdot {\bm u} & = 0
\end{aligned}
\]
where $\bm u = (v_1, v_2)$ represents the velocity, $p$ is the pressure, $\nu$ is the viscosity and $Re = 1/\nu$ is the Reynold's number. The domain of computation is $\Omega = [-0.5,1]\times[-0.5,1.5]$ with Direchilet boundary condition. We consider the Kovasznay flow problem that has the exact solution, i.e.,
\begin{align}
    v_1(x_1, x_2) &= 1 - e^{\lambda x_1} \cos(2 \pi x_2)\\
    v_2(x_1, x_2) &= \frac{\lambda}{2 \pi} e^{\lambda x_1} \sin(2 \pi x_2)\\
    p(x_1, x_2) &= \frac{1}{2} (1 - e^{2 \pi x_1} )
\end{align}
where $\lambda = \frac{1}{2\nu} - \sqrt{\frac{1}{4\nu^2} + 4\pi^2}$ and  the Reynold's number is set to 40. The Dirichlet boundary condition can be obtained by restricting the exact solution on the boundary of $\Omega$.

{\bf The Fokker-Planck equation} considered in case $(C_7)$ and $(C_8)$ is defined by
\begin{equation}\label{eq_adjoint}
\begin{aligned}
    \frac{\partial u(t,{\bm x})}{\partial t} + b(t, \bm x)\sum_{i=1}^{d}\frac{\partial u}{\partial x_{i}}(t,{\bm x}) + \frac{\sigma^2}{2}\sum_{i,j=1}^{d} \frac{\partial^2 u}{\partial x_{i}x_{j}}(t,{\bm x}) & = 0,\\
    u(0,\bm x) & = g(\bm x),
\end{aligned}
\end{equation}
where the coefficients $b(t,\bm x)$, $\sigma$, $g(\bm x)$ and the exact solutions are
\vspace{-0.2cm}
\begin{itemize}\itemsep0.1cm
    \item $(C_7)$: $b(x,t) = 2 \cos{(3 t)}$, $\sigma=0.3$, $u(x,0) = p(x; 0, 0.4^2)$ and $u(x,t) = p(x; \frac{2\sin{ (3t)}}{3}, 0.4^2 + t 0.3^2)$, where $p(x; \mu, \Sigma)$ denote the Gaussian density with mean $\mu$ and variance $\Sigma$.
    \item $(C_8)$: $b(x_1, x_2, t) = [\sin(2\pi t), \cos(2\pi t)]^T$, $\sigma=0.3$, $u(x_1, x_2, 0) = p(x; [0,0],  0.4^2 \bold{I}_2)$, and $u(x_1, x_2, t) = p(x; [-\frac{\cos(2\pi t)-1}{2\pi}, \frac{\sin(2\pi t)}{2\pi} ], (0.4^2 + t 0.3^2)\bold{I}_2)$, where $p(x; \mu, \Sigma)$ is the Gaussian density with mean $\mu$ and variance $\Sigma$.
\end{itemize}

{\bf The wave equation} considered in case $(C_9)$ is defined by
\[
\begin{aligned}
& \frac{\partial^2 u}{\partial t^2}  = c\frac{\partial^2 u}{\partial x^2},\;\; x \in [0,1], t\in [0,2] \\
& u(x,0) = \sin(4\pi x)\\
& u(0,t) = u(1,t)
\end{aligned}
\]
where $c = 1/(16\pi^2)$. The domain of computation is $\Omega = [0,1] \times [0,2]$; the exact solution is
\[
u(x,t) = \frac{1}{2} \left(\sin(4\pi x + t) + \sin(4\pi x - t)\right).
\]


\section{Setup of the experiments in Section \ref{sec:pde2}}\label{app:setup}
We specify the setup for the test cases $(C_1)$ to $(C_9)$ as follows:
\vspace{-0.2cm}
\begin{itemize}[leftmargin=20pt]\itemsep-0.0cm
    \item $(C_1)$: We evaluate the loss function in Eq.~\eqref{eq:ls} on a $50 \times 50$ uniform mesh in $\Omega = [-1,1]^2$, i.e., $J_1 = 2500$ in Eq.~\eqref{eq:ls}, and on 200 uniformly distributed points on $\partial \Omega$, i.e., $J_2 = 200$. After solving the least squares problem, we compute the error, i.e., the results shown in Figure \ref{fig:err} on a test set of 10,000 uniformly distributed random locations in $\Omega$. 
    \item $(C_2)$: We evaluate the loss function in Eq.~\eqref{eq:ls} on a $50 \times 50$ uniform mesh in $\Omega = [-1,1]^2$ and mask off the grid points outside the domain $\Omega = B_1(\bm 0)$, i.e., $J_1 = 1876$, and evaluate the boundary loss on 200 uniformly distributed points on $\partial \Omega$, i.e., $J_2 = 200$. After solving the least squares problem, we compute the error, i.e., the results shown in Figure \ref{fig:err} on a test set of 10,000 uniformly distributed random locations in $\Omega$.
    \item $(C_3)$: We evaluate the loss function in Eq.~\eqref{eq:ls} on a $50 \times 50$ uniform mesh in $\Omega = [-1,1]^2$ and mask off the grid points outside the domain $\Omega = [-1,1]^2 \backslash [0,1]^2$, i.e., $J_1 = 1875$, and evaluate the boundary loss on 200 uniformly distributed points on $\partial \Omega$, i.e., $J_2 = 200$. After solving the least squares problem, we compute the error, i.e., the results shown in Figure \ref{fig:err} on a test set of 10,000 uniformly distributed random locations in $\Omega$.
    \item $(C_4)$: We evaluate the loss function in Eq.~\eqref{eq:ls} on a $50 \times 50$ uniform mesh in $\Omega = [-1,1]^2$ and mask off the grid points outside the domain $\Omega = B_1(\bm 0) \backslash B_{0.5}(\bm 0)$, i.e., $J_1 = 1408$, and evaluate the boundary loss on 200 uniformly distributed points on $\partial \Omega$, i.e., $J_2 = 200$. After solving the least squares problem, we compute the error, i.e., the results shown in Figure \ref{fig:err} on a test set of 10,000 uniformly distributed random locations in $\Omega$.
    \item $(C_5)$: We evaluate the loss function in Eq.~\eqref{eq:ls} on a 10,000 uniformly distributed random locations in $\Omega = [-1,1]^3$,  i.e., $J_1 = 10000$, and evaluate the boundary loss on 2400 uniformly distributed points on $\partial \Omega$, i.e., $J_2 = 2400$, 400 points on each side of $\Omega$. After solving the least squares problem, we compute the error, i.e., the results shown in Figure \ref{fig:err} on a test set of 10,000 uniformly distributed random locations in $\Omega$.
    \item $(C_6)$: 
    We evaluate the loss function in Eq.~\eqref{eq:ls} on a $50 \times 50$ uniform mesh in $\Omega = [-0.5,1]\times[-0.5,1.5]$, i.e., $J_1 = 2500$ in Eq.~\eqref{eq:ls}, and on 200 uniformly distributed points on $\partial \Omega$ (50 points on each side of the box), i.e., $J_2 = 200$. We use Pichard iteration to handle the nonlinearity. Specifically, the residual loss is defined by
    \[
     loss = {\bm u_{\rm NN}^{k-1}} \cdot \nabla {\bm u_{\rm NN}^{k}} + \nabla p_{\rm NN}^k - \nu \Delta {\bm u_{\rm NN}^{k}}, 
    \]
    where $k$ is the Picard iteration number. In the $k$-th iteration, the nonlinear term ${\bm u_{\rm NN}^{k-1}} \cdot \nabla {\bm u_{\rm NN}^{k}}$ becomes linear due to the use of ${\bm u_{\rm NN}^{k-1}}$.
    After solving the least squares problem, we compute the error, i.e., the results shown in Figure \ref{fig:err} on a test set of 10,000 uniformly distributed random locations in $\Omega$. 
    \item $(C_7)$: The domain of computation is $(t,x) \in [0,1]\times[-2,2]$. We evaluate the loss function on a 50 (time) $\times$ 200 (space) = 10,000 grid points in the domain $\Omega$. We use the absorbing boundary condition in the spatial domain. We have a total of 3000 samples on the boundary of $\Omega$, i.e., 1000 samples for each of  $u(x,0)$, $u(2,t)$ and $u(-2,t)$. After solving the least squares problem, we compute the error, i.e., the results shown in Figure \ref{fig:err} on a test set of 10,000 uniformly distributed random locations in $\Omega$.
    \item $(C_8)$: The domain of computation is $t \in [0,1]$ and $(x_1,x_2) \in[-2,2]^2$. We evaluate the loss function on 10000 uniformly selected random points in the domain $\Omega$. We use the absorbing boundary condition in the spatial domain. In terms of samples on the boundary, we have $50 \times 50 = 2500$ grid points for the initial condition $u(x_1, x_2, 0)$, $20(\text{time}) \times 50(\text{space}) = 1000$ grid points for each of $u(\pm 2, x_2, t)$ and $u(x_1,\pm 2, t)$. After solving the least squares problem, we compute the error, i.e., the results shown in Figure \ref{fig:err} on a test set of 10,000 uniformly distributed random locations in $\Omega$.
    \item $(C_9)$: We evaluate the loss function in Eq.~\eqref{eq:ls} on $50\text{(time)} \times 100\text{(space)} = 2500$ grid points in domain,  i.e., $J_1 = 10000$, and evaluate the boundary loss on 1000 uniformly distributed points on $\partial \Omega$, i.e., $J_2 = 1500$, 500 points on each side of $\Omega$. After solving the least squares problem, we compute the error, i.e., the results shown in Figure \ref{fig:err} on a test set of 10,000 uniformly distributed random locations in $\Omega$.
\end{itemize}

We use the standard least squares solver \texttt{torch.linalg.lstsq} in Pytorch to solve all the least squares problems. Our code is implemented using Pytorch on a workstation with an NVIDIA Tesla V100 GPU. 

{\bf Setup for PINN}. 
For each test case, PINN uses exactly the same setting as TransNet, including network architecture, loss function, data, to ensure fair comparison. In terms of training, we set learning rate to 0.001 with a decrease factor of 0.7 every 1000 epochs. We first use Adam optimizer to train the neural networks for 5000 epochs, which gives us the results in Figure \ref{fig:err} labeled by ``PINN:Adam''. Then we continue training the network using LBFGS for another 200 iterations, which gives us the results in 
Figure \ref{fig:err} labeled by ``PINN:Adam+BFGS''.

{\bf Setup for for the random feature models}. The random feature model use 
exactly the same setting as TransNet, including network architecture, loss function, data, to ensure fair comparison. The parameters $\{\bm w_m, b_m\}_{m=1}^M$ are determined by the default initialization methods in Pytorch, and the parameters in the output layer is obtained by the least squares solver \texttt{torch.linalg.lstsq} in Pytorch.

\end{document}